\documentclass[a4paper,11pt]{article}

\usepackage{pstricks}
\usepackage{amsfonts,latexsym,amsmath,amssymb,graphicx,url}

\def\kaxxa{{\vcenter {\hrule height .2mm
\hbox{\vrule width .2mm height 2mm \kern 2mm
\vrule width .2mm} \hrule height .2mm}}}
\def\lesta{\hfill $\kaxxa$ \medskip}

\newcommand{\diam}{\mbox{diam}}
\newcommand{\rn}{\mbox{rn}}
\newcommand{\Crn}{{\cal C}\mbox{rn}}
\newcommand{\ern}{\mbox{ern}}
\newcommand{\Cern}{{\cal C}\mbox{ern}}
\newcommand{\Tern}{{\cal T}\mbox{ern}}
\newcommand{\Trn}{{\cal T}\mbox{rn}}

\newcommand{\ED}{{\cal{ED}}}
\newcommand{\D}{{\cal {D}}}

\newtheorem{Lemma}{Lemma}[section]

\newtheorem{Theorem}{Theorem}[section]

\newtheorem{Observation}{Observation}[section]
\newtheorem{Conjecture}{Conjecture}[section]

\title{On the edge-reconstruction number of a tree}
\author{K. Asciak\footnote{kevin.j.asciak@um.edu.mt} and J. Lauri\footnote{josef.lauri@um.edu.mt} \\Department of Mathematics \\ University of Malta \\ Malta \and W.~Myrvold\footnote{wendym@cs.uvic.ca} \\ Dept. of Computer Science\\
University of Victoria\\ Victoria, B.C.\\ Canada V8N 6K3
\and V. Pannone\footnote{virgilio.pannone@unifi.it} \\Dipartment of Mathematics and Informatics\\ University of Florence\\ Italy.}

\begin{document}

\maketitle

\begin{abstract}
The edge-reconstruction number $\ern(G)$ of a graph $G$ is equal to the minimum number of edge-deleted subgraphs $G-e$ of $G$ which are sufficient to determine $G$ up to isomorphsim. Building upon the work of Molina and using results from computer searches by Rivshin and more recent ones which we carried out, we show that, apart from three known exceptions, all bicentroidal trees have edge-reconstruction number equal to 2. We also exhibit the known trees having edge-reconstruction number equal to 3 and we conjecture that the three infinite families of unicentroidal trees which we have found to have edge-reconstruction number equal to 3 are the only ones.
\end{abstract}

\section{Introduction}

Trees have often been the test-bed for various graph theoretic conjectures, not least being the Reconstruction Conjecture. Kelly's proof that trees are reconstructible \cite{kelly57} was the first substantial reconstructibility proof. This result was later improved by various authors who showed that trees can be reconstructed using only their endvertex- or peripheral-vertex- or cutvertex-deleted subgraphs \cite{harary&pal66,bondy69,lauri83}.

A vertex-deleted subgraph $G-v$ of $G$ is called a \emph{card} of $G$; the collection of cards of $G$ is called the \emph{deck} of $G$, denoted by $\D(G)$. Our main focus in this paper will be on the analogously defined \emph{edge-cards} of $G$ which are the edge-deleted subgraphs $G-e$ of $G$; the collection of edge-cards of $G$ is called the \emph{edge-deck} of $G$ and is denoted by $\ED(G)$.

In \cite{harary&pla85}, Harary and Plantholt introduced the notion of reconstruction numbers. The \emph{reconstruction number} $\rn(G)$ of a graph $G$ is defined to be the least number of vertex-deleted subgraphs of $G$ which alone reconstruct $G$ uniquely (up to isomorphism). The \emph{class reconstruction number} $\Crn(G)$ is defined as follows. Let $\cal C$ be a class of graphs closed under isomorphism. Then the class reconstruction number of a graph $G$ in $\cal C$ is the minimum number of vertex-deleted subgraphs of $G$ which, together with the information that $G$ is in $\cal C$, reconstruct $G$ uniquely.
It is clear that the reconstruction number of a graph is always at least 3 and that $\Crn(G)\leq\rn(G)$. In fact, the class reconstruction number can even be 1, for example, when $\cal C$ is the class of regular graphs.
The \emph{edge-reconstruction number} $\ern(G)$ of a graph $G$ and the \emph{class edge-reconstruction number} $\Cern(G)$ for a graph $G$ in $\cal C$ are analogously defined.

In \cite{harary&lau88}, Harary and Lauri tackled the reconstruction number of a tree. Let $\cal T$ be the class of trees. In their paper, Harary and Lauri tried to show that $\Tern(T)\leq2$. Although they managed to achieve this in many of the cases they considered, in some cases they had to settle for the upper bound of 3. So, what was accomplished in \cite{harary&lau88} was to show that $\Trn(T)\leq3$ and to make plausible their conjecture that, in fact, $\Trn(T)\leq2$ for all trees $T$. Myrvold \cite{myrvold90} soon improved the first result by showing that
$rn(T)\leq3$. The conjecture $\Trn(T)\leq2$, however, still stood. A significant step forward was recently taken by Welhan \cite{welhan10} who proved that the class reconstruction number of trees is at most 2 for trees without vertices of degree 2.

The situation for the edge-reconstruction numbers of trees is less clear, somewhat surprisingly compared with what happens in the Reconstruction Problem where edge-reconstruction is easier than vertex-reconstruction. Although Harary and Lauri conjectured that $\Trn(T)\leq2$ for all trees $T$, they presented in \cite{harary&lau88} a few trees with class edge-reconstruction number $\Tern$ equal to 3 even though their class (vertex) reconstruction number was equal to 2. In \cite{molina93}, Molina started to tackle the edge-reconstruction number of trees. In summary, these are Molina's main results.

\begin{enumerate}
\item Let $T$ be a unicentroidal tree with at least four edges, then $\ern(T) \leq 3$.

\item Let $T$ be bicentroidal with centroidal vertices $a$ and $b$, and let $G$ and $H$ be the two components of $T-ab$ with $a$ in $G$ and $b$ in $H$. Then

\begin{enumerate}
 \item If one of the centroidal vertices has degree equal to two, then $\ern(T) \leq 3$.
 \item If both centroidal vertices have degree at least three and if $G$ or $H$ has an irreplaceable endvertex (defined below), then $\ern(T)= 2$.
 \item If both centroidal vertices have degree at least three and if either $G$ or $H$ has no irreplaceable endvertex, then $\ern(T) \leq 3$.
\end{enumerate}
\end{enumerate}

In this paper we shall improve the above results on bicentroidal trees by showing that $\ern(T)=2$ when the degrees of the centroidal vertices are 2 and even when both $G$ and $H$ have no irreplaceable vertices, giving our main result is that all bicentroidal trees, with only three exceptions, have $\ern$ equal to 2. We shall also prove some results on unicentroidal trees and, based on these results and empirical evidence which we shall present, we give a conjecture stating which infinite classes of unicentroidal trees have $\ern$ equal to 3.

%\begin{eqnarray}
%$3 & = & 1+2\\
%7 & = & 3+4$
%
%\end{eqnarray}

One final definition: suppose we are considering $\rn(G)$ or $\ern(G)$ and suppose that a graph $H\not\simeq G$ has in its deck (edge-deck) the cards (edge-cards) $G-v_1,\ldots,G-v_k$ ($G-e_1,\ldots,G-e_k$) we then says that $H$ is a \emph{blocker} for these cards (edge-cards) or that $H$ \emph{blocks} these cards (edge-cards).

\section{Main techniques}

We shall here present the main techniques and supporting results used in this paper. Many of these were first used or proved in \cite{harary&lau88}. While all work on the reconstruction of trees prior to \cite{harary&lau88} depended on the centre of a tree, in \cite{harary&lau88} the centroid was used instead. Since then, all investigations of reconstruction numbers of trees depended heavily on centroids. Non-pseudosimilarity and irreplaceabilty of endvertices were also very important techniques first used in the proofs in \cite{harary&lau88}. These ideas will be explained below. We shall also present a new technique and a result which will be used for the first time in this paper.

\subsection{The centre and the centroid of a tree, rooted trees and branches}

The \emph{diameter} $\diam(G)$ of a connected graph $G$ is the length of a longest path in $G$. The \emph{eccentricity} of a vertex $v$ in $G$ is the longest distance from $v$ to any other vertex in the graph. The \emph{centre} of $G$ is the set of vertices with minimum eccentricity. It is well-known that if $G$ is a tree then the centre is either one vertex or two adjacent vertices. 

\medskip\noindent
We now turn our attention to the centroid. Define the \emph{weight} of a vertex $v$ of a tree $T$, denoted by $wt(v)$, to be the number of vertices in a largest component of $T - v$. For example all endvertices in a $n$-vertex tree have weight $n-1$. The \emph{centroid} of a tree $T$ is the set of all vertices with minimum weight denoted by $wt(T)$. A \emph{centroidal vertex} is a vertex in the centroid. It is well-known that the centroid of a tree consists of either one vertex or two adjacent vertices. A tree with one centroidal vertex is called \emph{unicentroidal} while a tree with two centroidal vertices is called 
\emph{bicentroidal}. In the latter case, the edge joining the centroidal vertices is called the \emph{centroidal edge}. When $T$ is bicentroidal with centroidal edge $e$, the two components of $T-e$ are also said to be \emph{centroidal components}.

The following simple observation will be very useful. The second part, especially, tells us that for a graph $T$ which we know to be a tree, if it is bicentroidal, then one can determine from an edge-deleted subgraph $T-e$ of $T$ alone, whether or not $e$ is the centroidal edge of $T$ and also, if $e$ is the centroidal edge, the isomorphism type of the two centroidal components.

\begin{Observation} \label{obs:recogbicentroidal}
Let $T$ be a tree of order $n$ and let $v$ be a vertex of $T$. Then $wt(v)\leq \frac{n}{2}$ if and only if $v$ is in the centroid of $T$. Also, $T$ is bicentroidal with centroidal vertices $a$ and $b$ if and only if $T - ab$ has two components $G,H$ each of order $\frac{\left|V(T)\right|}{2}$.
\end{Observation}

\noindent
{\bf Notation.} In the rest of the paper, $a$ and $b$ will denote the centroidal vertices of a bicentroidal tree with centroidal components $G$ and $H$ such that $a$ is in $G$ and $b$ is in $H$.

%Let $T$ be a bicentroidal tree whose centroidal vertices are $a$ and $b$. Then $T_1$ is the %largest component of $T - c_2$ while $T_2$ is the largest component of $ T - c_1$. Clearly, %$c_1 \in V(T_1)$ and $c_2 \in V(T_2)$. 

\medskip\noindent
A \emph{rooted tree} is a tree which has one identified vertex. Let $P$ be a path in a tree and let $v$ be an internal vertex on $P$. The \emph{branch} at $v$ relative to $P$ is the subtree, rooted at $v$, induced by all those vertices connected to $v$ by a path not containing other vertices of $P$. 

A vertex of degree 1 is said to be an \emph{endvertex}. A cutvertex in a tree which is adjacent to only one vertex of degree greater than 1 is said to be an \emph{end-cutvertex}. An edge incident to an endvertex is called an \emph{end-edge}.

\subsection{Pseudosimilar vertices, irreplaceable edges and conjugate pairs of trees}

Most of the works which we mentioned and which deal with reconstruction numbers of trees of some sort make heavy use of the impossibility of endvertices being pseudosimilar in a tree and of the fact that only a few very special type of trees have the property that any end-edge can be exchanged with another giving us a tree isomorphic to the one which we started with. Since we shall be using these results even in this paper we shall explain them and their general use in this section. We shall also prove another result in this vein which we shall be needing, namely a result about a pair of trees such that any one can be obtained from the other by exchanging some end-edges in a particular way

Let $u$ and $v$ be two vertices in a graph $K$ such that an automorphism of $K$ maps $u$ into $v$. Then $u$ and $v$ are said to be \emph{similar} in $K$. Now suppose that $u$ and $v$ are such that $K - u$ is isomorphic to $K - v$; we call such a pair of vertices \emph{removal-similar}. If $u$ and $v$ are removal-similar in $K$ but not similar, then $u$ and $v$ are said to be \emph{pseudosimilar vertices} in $K$. The following results  say that endvertices and end-cutvertices in a tree cannot be pseudosimilar.

\begin{Theorem}\textbf{(Harary and Palmer)} \cite{harary&pal66a} (i) Any two removal-similar endvertices in a tree are similar.\\
\textbf{(Kirkpatrick, Klawe and Corneil)} \cite{kirkpatrick&83} (ii) Any two removal-similar end-cutvertices in a tree are similar. 
\end{Theorem}

Since we shall be expanding on this and the subsequent result in this paper it is interesting to see one way in which these two results have been extended by Krasikov in \cite{krasikov91}. Let $T$ be a tree and $a,b\in V(T)$, and let $A,B$ be two rooted trees. Then $T_{a,b}(A,B)$ denotes the tree obtained by identifying the root of $A$ with $a$ and the root of $B$ with $b$. Krasikov proved the following.

\begin{Theorem}
If $A$ and $B$ are two non-isomorphic rooted trees and 
\[T_{a,b}(A,B) \simeq T_{a,b}(B,A)\] 
then $a$ and $b$ are similar in $T$.
\end{Theorem}

Clearly, if we take $A$ to be the tree on two vertices and $B$ a single vertex, then this result gives that endvertices cannot be pseudosimilar in a tree.

\begin{figure} 
 \centering
 \includegraphics[width=8.5cm, height=4cm]{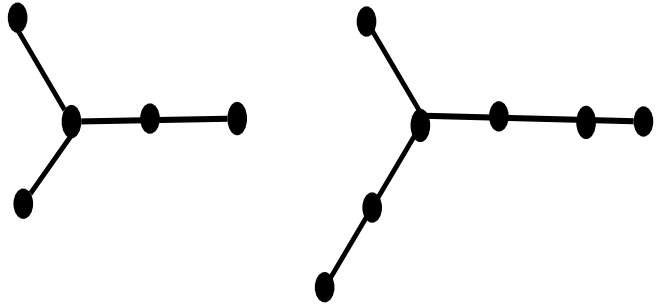}
 \caption{The trees: (a) $S_1$; and (b) $S_2$} \label{fig:pseudopaths}
 \end{figure}

\medskip\noindent
Now let $e=xv$ be an end-edge of $T$ with $\deg(v)=1$. Let $y\not=x$ be another vertex of $T$ and let $T'= T - e + e'$, where $e'=yv$. If $T'$ is isomorphic to $T$, then $e$ is called a \emph{replaceable end-edge}. If there is no such vertex $y$ then $e$ is called an \emph{irreplaceable end-edge}. Let $S_1$ and $S_2$ be the graphs shown in Figure \ref{fig:pseudopaths}. A tree which is isomorphic either to a path $P_k$ on $k$ vertices or to one of $S_1$ or $S_2$ is said to be a \emph{pseudopath}.

The following theorem was proved in \cite{harary&lau88} and was also profitably used in \cite{molina93}.

\begin{Theorem} \label{thm:replaceable} Any tree which is not a pseudopath has an irreplaceable end-edge.
\end{Theorem}

The use of non-pseudosimilarity of endvertices and irreplaceable edges are important techniques 
which are used in these two broad scenarios in this paper. First of all, suppose that we have two trees $G, H$ and we know that the tree $T$ to be reconstructed is obtained by joining together with a new edge an endvertex $a$ of $G$ to another endvertex $b$ of $H$ (we do not know which vertices are $a$ and $b$). Suppose, however, that we know the isomorphism types of both $G'=G-a$ and 
$H'=H-b$. Then, since endvertices in a tree cannot be pseudosimilar, we can pick any endvertex $x$ in $G$ such that $G-x \simeq G'$ and similarly any endvertex $y$ in $H$ such that $H-y\simeq H'$, and join the two vertices $x$ and $y$ giving the reconstruction of $T$ which is unique up to isomorphism. 

The second scenario is basically this. Suppose that we know again that the tree $T$ to be reconstructed is obtained by joining vertex $a$ in $G$ to vertex $b$ in $H$ ($a, b$ need not be endvertices now). We are also given the tree $T'$ which is composed of $G$ joined correctly to $H'$, where $H'$ is $H$ less an endvertex and we can identify the edge $ab$ in $T'$. We therefore know from $T'$ the components $G$ and $H'$ and how they are connected. We just need to be able to put back the missing endvertex in $H'$. In order to have unique reconstruction up to isomorphism, non-pseudosimilarity of the missing endvertex is not enough here. We now require that the missing vertex be irreplaceable in $H$.

\medskip\noindent
In this paper we shall also need a notion which is in some way an extension of the idea of replaceable endvertices. Instead of asking that exchanging an end-edge in a tree gives us the same tree, we ask that a pair of trees are related by a particular exchange of end-edges. This is quite a natural occurrence when considering reconstruction of trees. First we need a technical definition which, however, will find its natural place in our reconstruction results later in Theorem \ref{thm:main2}. 

Suppose $G$ and $H$ are two \emph{non-isomorphic trees}. Let $a, b$ be endvertices of $G$ and $H$, respectively. Suppose also that:
\begin{enumerate}
\item $G-a + e_1 \simeq H$ for some new end-edge $e_1$ added to $G-a$;
\item $G+aa' - e_2 \simeq H$ for some new endvertex $a'$ added to $G$ and some end-edge $e_2$ of $G$;
\item $H-b + e_3 \simeq G$ for some new end-edge $e_3$ added to $H-b$;
\item $H+bb' - e_4 \simeq G$ for some new endvertex $b'$ added to $H$ and some end-edge $e_4$ of $H$.  
\end{enumerate} 
Then $G$ and $H$ are said to be a \emph{conjugate pair} of trees.

\bigskip\noindent
The theorem we shall need is the following.

\begin{figure} 
 \centering
 \includegraphics[width=8.5cm, height=4cm]{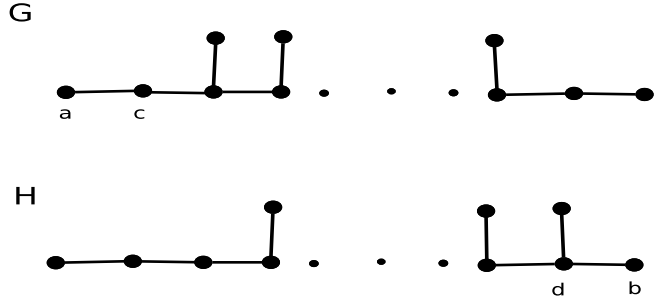}
 \caption{A conjugate pair of trees} \label{fig:pannone}
 \end{figure}

\begin{Theorem} \label{thm:pannone}
Let $G$ and $H$ be a conjugate pair of trees as in the definition. Then $G$ and $H$ must be trees as shown in Figure \ref{fig:pannone}.
\end{Theorem}

\noindent
{\bf Proof.} Let $c$ be the neighbour of $a$ in $G$ and $d$ the neighbour of $b$ in $H$. We shall consider two cases:

\medskip\noindent \emph{Case 1: At least one of the trees on the left-hand-side of equations (1)--(4) in the definition of conjugate pairs has a different centre from the original tree on the right-hand-side.}

\smallskip\noindent
We shall suppose that the change of centre occurs in Equation 1 of the definition. The arguments for Equation 2 are similar, and those for Equations 3 and 4 follow by symmetry. Therefore we have that the centre of $G-a + e_1$ is not the same as in $G$. We shall consider the case when the centre of $G$ has one vertex. The bicentral case is similar.

\begin{figure} 
 \centering
 \includegraphics[width=9.5cm, height=3.5cm]{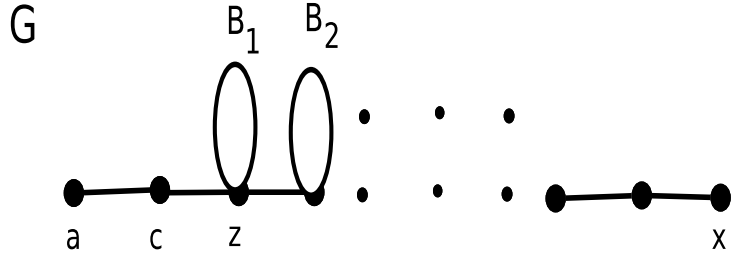}
 \caption{The tree $G$ and the branches relative to $P$} \label{fig:pannone2}
 \end{figure}

The condition implies that $a$ is on the unique longest path of $G$. Call this path $P$. Let $x$ be the vertex at the other end of $P$. It also follows that $\deg(c)=2$. If $z$ is the other neighbour of $c$ on $P$, then the branch at $z$ relative to $P$ must be at most an edge; the branch at the next vertex $z'$ cannot have a vertex distant more than 2 from $z'$, and so on. Let the branches relative to $P$ at the internal vertices of $P$ be $B_1, B_2, \ldots, B_k$ in that order starting from the branch at $z$ as shown in Figure \ref{fig:pannone2}. Since, by definition, $G$ and $H$ are not isomorphic, the branches cannot all be trivial (consisting of only the root vertex).

Also, from Equation 1 it follows that $\diam(H)\leq \diam(G)$ and, from Equation 2, that $\diam(G)\leq \diam(H)$. Therefore $G$ and $H$ have the same diameter. Therefore the end-edge $e_1$ in Equation 1 of the definition is $xx'$ for some new vertex $x'$. 

\begin{figure} 
 \centering
 \includegraphics[width=9.5cm, height=3.5cm]{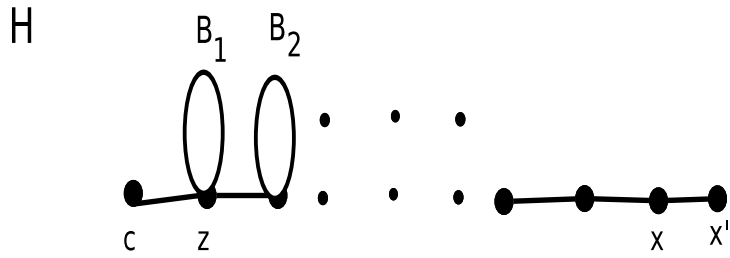}
 \caption{The tree $H$ shown as $G-a+xx'$} \label{fig:pannone3v2}
 \end{figure}

Now, consider $H$ given as $G-a + xx'$ as depicted in Figure \ref{fig:pannone3v2}. Which would be the vertex $b$ in $H$ which satisfies Equation 4? Recall that $H+bb' - e_4$, being isomorphic to $G$, would have to have a unique path of maximum length and the branches of the internal vertices of this path would have to be $B_1, B_2,\ldots, B_k$ in that order.

Therefore $b$ cannot be $x'$ nor any endvertex in any of $B_1, B_2, \ldots, B_k$. Therefore $b$ must be the vertex $c$ and $d$ must be the vertex $z$ in Figure \ref{fig:pannone3v2}.

\begin{figure} 
 \centering
 \includegraphics[width=9.5cm, height=3.5cm]{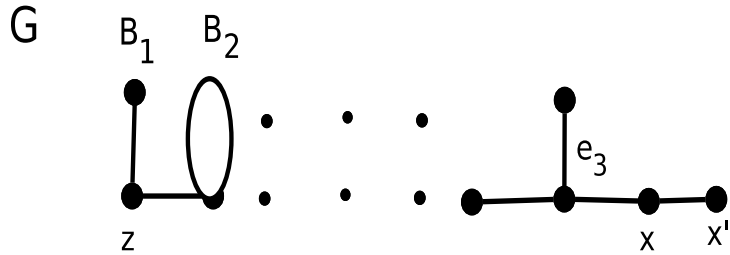}
 \caption{The graph $G$ as $H-b+e_3 \simeq H-c + e_3$} \label{fig:pannone4v2}
 \end{figure}

But, in order to satisfy Equation 3, the branch $B_1$ must be a single edge and the end-edge $e_3$ must be as shown in Figure \ref{fig:pannone4v2}.   

But then, comparing $G$ as in Figure \ref{fig:pannone4v2} with Figure \ref{fig:pannone2} shows that all the branches are single edges and $G$ is as in Figure \ref{fig:pannone}. This finishes Case 1.

\bigskip\noindent \emph{Case 2: Every tree on the left-hand-side of equations (1)--(4) in the definition of conjugate pairs has the same centre as that in the original tree on the right-hand-side.}

\medskip\noindent
We shall prove that this leads to a contradiction, therefore only Case 1 can hold. 
We shall only consider the unicentral case. The bicentral case can be treated similarly. First we need to define exactly what we mean by a central branch of a central tree $T$ with central vertex $v$. Let $A$ be a component of $T-v$ and let $u$ be the neighbour of $v$ in $A$. Let $B=A+uv$. Then $B$ will be called a \emph{central branch} of $T$. Clearly, the number of central branches of $T$ is equal to $\deg(v)$.

Consider first Equation 1: $G-a + e_1 \simeq H$. Let $B_1$ be the central branch of $G$ containing the edge $ac$. We now have two sub-cases.

\bigskip\noindent{Case 2.1: The edge $e_1$ is incident to a vertex in $B_1$.}

\medskip\noindent
Therefore $G$ and $H$ have exactly the same collection of branches except that $H$ has the branch $B_1' \simeq B_1 -a + e_1$ instead of $B_1$. So, the only way of obtaning $G$ back from $H$ in the way stipulated by Equations 3 and 4 is by changing $B_1'$ to $B_1$ and, similarly, the only way of going from $G$ to $H$, according to Equations 1 and 2, is by changing $B_1$ to $B_1'$. Therefore Equations 1 to 4 hold for the trees $B_1$ and $B'_1$, that is, they form a conjugate pair of trees. Applying induction on the number of vertices gives us that $B_1$ and $B'_1$ are as specified by the theorem, that is, as in Figure \ref{fig:pannone}. Therefore $G$ is the tree $B_1$ with extra branches joined to $v$ (which is an endvertex in $B_1$). But then, $G$ cannot satisfy Equations 1 to 4, that is, it cannot be a member of a conjugate pair of trees.

\bigskip\noindent{Case 2.2: The edge $e_1$ is not incident to a vertex in $B_1$.}

\medskip\noindent
Let $B_2$ be the central branch of $G$ containing $e_1$. Therefore the branches of $G$ and $H$ are identical except that $H$ has $B_1'=B_1-c$ instead of $B_1$ and $B_2' = B_2+e_1$ instead of $B_2$.

Now, the endvertex $a$ of $G$ which is in $B_1$ is also involved in Equation 2: $G+aa'-e_2 \simeq H$. Let us consider where the edge $e_1$ can come from so that $H$ is isomorphic to both $G-a+e_1$ and $G+aa'-e_2$. We point out that we need to obtain the same collection of branches for $H$ (with $B_1'$ and $B_2'$ instead of $B_1$ and $B_2$, respectively) and that we cannot do this by moving the centre. That is, we can only make modifications to the existing central branches. 

The only way this can happen is if $e_2$ comes from some third central branch $B_3$. Now consider the orders of $B_1, B_2, B_3$. Let these orders be $r,s,t$, respectively. Then, a moment's consideration shows that we must have that $r=p+1, s=p$ and $t=p+2$, for some $p$. 

\begin{figure} 
 \centering
 \includegraphics[width=8.5cm, height=4cm]{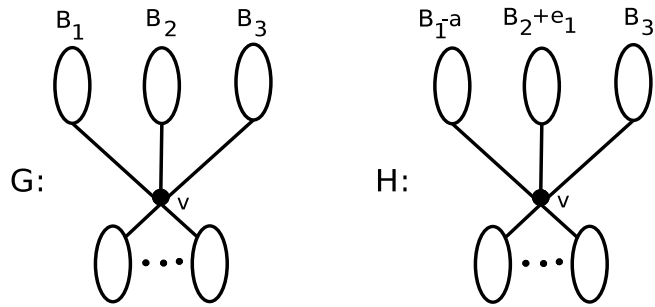}
 \caption{The graphs $G$ and $H=G-a + e_1$} \label{fig:pannone5}
 \end{figure}

Therefore $G$ and $H$ are as shown in Figure \ref{fig:pannone5}. 

\begin{figure} 
 \centering
 \includegraphics[width=8.5cm, height=5cm]{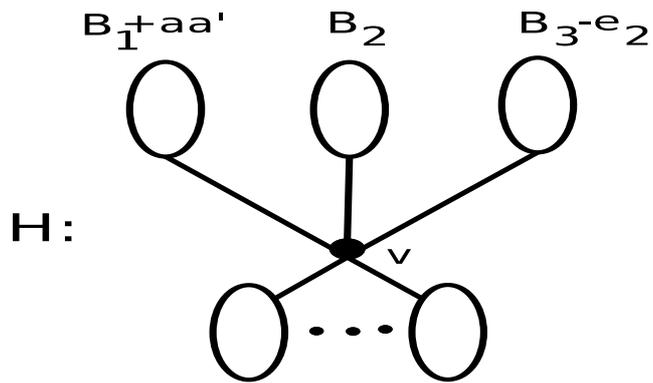}
 \caption{The graph $H=G+aa'-e_2$} \label{fig:pannone6}
 \end{figure}

But also, $H$ is isomorphic to the tree shown in Figure \ref{fig:pannone6}. So we get, for example, by considering orders, that $B_3\simeq B_1+aa'$ and $B_2\simeq B_1-a$. Switching over from $G$ to $H$ and from $H$ to $G$ using Equations 1 to 4 involves exchanges endvertices between these three branches (or three branches in $G$ or $H$ isomorphic to them). 

\begin{figure} 
 \centering
 \includegraphics[width=8.5cm, height=4cm]{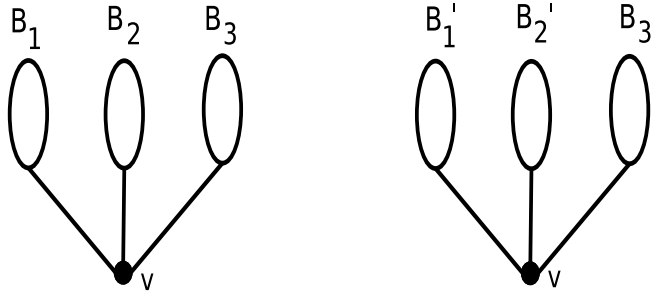}
 \caption{These two subtrees must be isomorphic} \label{fig:pannone7}
 \end{figure}

So, when we are considering $G$ in Equations 1 and 2, the vertex equivalent to $a$ would be in that branch which has order $p+1$, the new edge $e_1$ would be attached to the branch of order $p$, and $e_2$ would be removed from the branch of order $p+2$. Similarly, if we are considering $H$, for $b, e_3$ and $e_4$ in Equations 3 and 4. But we are always permuting between the same (up to isomorphism) three branches which become isomorphic to $B_1, B_2, B_3$ in $G$ and $B_1', B_2', B_3$ in $H$. But this would force the two trees in Figure \ref{fig:pannone7} to be isomorphic, therefore $G$ and $H$ would be isomorphic, a contradiction which completes our proof.
\lesta

\begin{figure} 
 \centering
 \includegraphics[width=8.5cm, height=4cm]{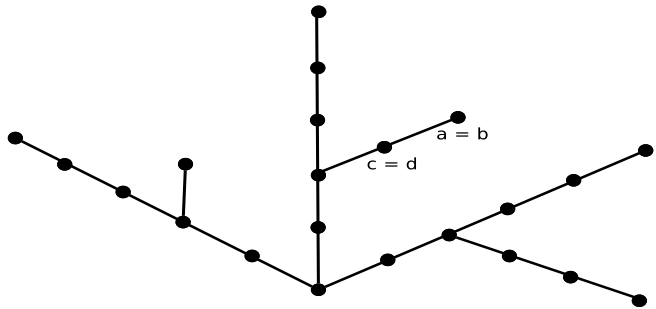}
 \caption{$G\simeq H$ would be a conjugate pair if allowed to be isomorphic}  \label{fig:pannone8}
 \end{figure}

\noindent
Note that it is the fact that $G$ and $H$ are not isomorphic which forces conjugate pairs to be as described in the theorem and which gives us our final contradiction. If $G$ and $H$ are allowed to be isomorphic then, for example, two trees both isomorphic to the one shown in Figure \ref{fig:pannone8} do satisfy Equations 1 to 4. Note, in this example, the three central branches as described in Case 2.2 of the above proof.

\subsection{Recognising trees and Molina's Lemma}

There is a simple but very useful result proved by Molina in \cite{molina93} which often allows us to identify a graph as a tree from two given edge-cards. We reproduce its short proof for completeness' sake.

\begin{Lemma} \label{lem:molina}
Let $G$ be a  graph with edges $e_1$ and $e_2$. Suppose that the edge-card $G-e_1$ has two  components which are trees of orders $p_1$ and $p_2$ while the edge-card $G-e_2$ has another two components which are trees of orders $q_1$ and $q_2$. If $\{p_1,p_2\}\not =\{q_1,q_2\}$, then $G$ is a tree.
\end{Lemma}

\noindent
{\bf Proof.} Suppose $G$ is not a tree. 
Without loss of generality we assume that $e_1$ joins two vertices in the same component of $G-e_1$; call this component $H$, that is, $H+e_1$ contains a cycle. Therefore to obtain the second 
edge-card with two trees as components an edge must be removed from $H+e_1$. But this contradicts that $\{p_1,p_2\}\not=\{q_1,q_2\}$. \lesta

\subsection{Some special types of tree}

\begin{figure} 
 \centering
 \includegraphics[width=8.5cm, height=4cm]{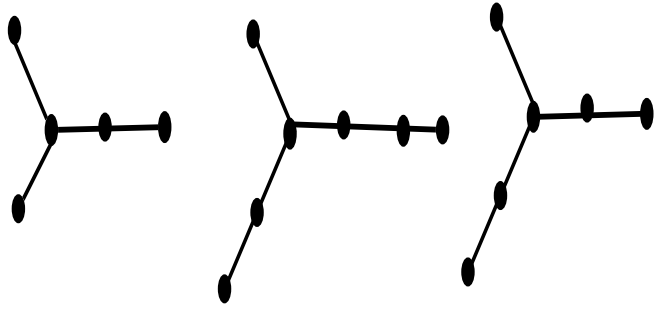}
 \caption{The trees: (a) $S_{1,1,2} (=S_1)$; (b) $S_{1,2,3} (=S_2)$; and (c) $S_{1,2,2}$} \label{fig:pseudostars}
 \end{figure}

A special type of tree denoted by $S_{p,q,r}$ is a unicentroidal tree similar to a star (that is, the tree on $n$ vertices, $n-1$ of which are endvertices) which consists of three paths on $p$, $q$ and $r$ edges, respectively, emerging from the centroidal vertex. Some examples are shown in Figure \ref{fig:pseudostars}. Note that the pseudopaths $S_1$ and $S_2$ defined above are $S_{1,1,2}$ and $S_{1,2,3}$, respectively.

\begin{figure} 
 \centering
 \includegraphics[width=9cm, height=2cm]{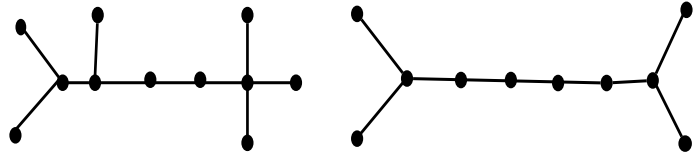}
 \caption{Caterpillars: (a) $C(2,1,0^2,3)$; (b) $C(2,0^4,2)$} \label{fig:caterpillars}
 \end{figure}

A \emph{caterpillar} is a tree such that the removal of all of its endvertices results in a path.
This path is called the \emph{spine} of the caterpillar. A caterpillar whose spine is the
path $v_1v_2\ldots v_s$ and such that the vertex $v_i$ is adjacent to $a_i$ endvertices will be denoted by 
$C(a_1, ..., a_s)$. Two examples are shown in Figure \ref{fig:caterpillars}. Finally, a path on $n$ vertices is denoted by $P_n$.

\section{Bicentroidal trees}

\subsection{The centroidal component $G$ is not a pseudopath and $\deg(b)\geq3$}

\begin{Theorem} \label{thm:main1}
Let $T$ be a bicentroidal tree with bicentroidal edge $ab$, bicentroidal components $G,H$, $a\in V(G)$, $b\in V(H)$. Suppose $\deg(b)\geq3$ and $G$ is not a pseudopath.
Then $\ern(T)=2$.
\end{Theorem}

\noindent
{\bf Proof.} Since $T$ is bicentroidal and $ab$ is the centroidal edge then the two components $G$ and $H$ of the card $T - ab$ have the same number of vertices, namely $\frac{\left|V(T)\right|}{2}$.
Let $f$ be an irreplaceable end-edge of $G$ (such an $f$ exists since $G$ is not a pseudopath). We claim that $T$ is reconstructible from $T-ab$ and $T-f$.

By Lemma \ref{lem:molina} we can recognise from $T-ab$ and $T-f$ that the graph to be reconstructed is a tree. 
By Observation \ref{obs:recogbicentroidal} one can therefore recognise from the edge-card $T - ab$ that the edge $ab$ is the centroidal edge and also that $G$ and $H$ are the centroidal components of $T$.

Now, we would like to show
that the centroidal edge is recognisable in the edge-card $T - f$.  There is surely an edge $e$ such that $(T - f) - e$ has non-trivial components $G - f$ and $H$, but we can definitely say that $e$ is the edge $ab$ only if:
\begin{itemize}
\item[(i)] there is only one edge $e$ such that the non-trivial components of $(T-f)-e$ are
isomorphic to $H$ and some $T-f$; 

and

\item[(ii)]
there is no edge $e'$ such that the non-trivial component of $(T - f) - e'$ is isomorphic to 
$G$ and some $H - f$.
\end{itemize}

If both (i) and (ii) hold then we can distinguish the centroidal edge in $T - f$ and we can reconstruct uniquely by putting $f$ back into $G - f$, since $f$ is an irreplaceable end-edge (note that this proof also works if the end-edge $f$ happens to be adjacent to the centroidal edge). 

But cases (i) and (ii) can fail to occur only if the degree of the centroidal vertex $b$ is two. Since $\deg(b)>2$ it follows that $T$ is reconstructible from $T-ab$ and $T-f$. \lesta

We shall come back to what happens when $\deg(b)=2$ but $G$ is still not a pseudopath in Lemma \ref{lem:remains} after having obtained some more results and discussed some special cases.

\subsection{Both $\deg(a)$ and $\deg(b)$ equal 2 and none of $G$ or $H$ is a pseudopath}

\begin{Theorem} \label{thm:main2}
Let $T$ be a bicentroidal tree with bicentroidal edge $ab$. Let $\deg(a)=\deg(b)=2$ and suppose that none of the two centroidal components $G$, $H$ is a pseudopath. Then $\ern(T)=2$. 
\end{Theorem} 
 
\noindent
{\bf Proof.} Recall that $b\in V(H)$; let $d\in V(H)$ be the other neighbour of $b$.  We shall first try to show that $T$ is reconstructible from $T-ab$ and $T-bd$ and we shall see where this can go wrong.

As before, from the two given edge-cards we can recognise that $T$ is a tree and that $G,H$ are its centroidal components. Consider $T-bd$. If we can definitely tell that the larger component of $T-bd$ is $G$ plus some edge then we would only need to decide which is the extra end-edge in the larger component. But, since endvertices cannot be pseudosimilar, we can choose any endvertex whose deletion gives $G$. We therefore know, up to isomorphism, which of the vertices of $G$ is incident to the centroidal edge. Now we would need to do the same with $H$. 

Let $H'$ be the smaller component of $T-bd$. Recall that we know the component $H$. We look for any endvertex $d'$ such that $H-d' \simeq H'$. Again, by non-pseudosimilarity of endvertices, any such choice is equivalent to $d$ up to isomorphism. So we also know the vertex of $H$ which is incident to the centroidal edge, hence $T$ can be uniquely reconstructed. 

This proof fails if we cannot tell whether the larger component is $G$ plus an end-edge or $H$ plus an end-edge. This ambiguity can only happen if $G\not\simeq H$ and $G+ab-\alpha \simeq H$ for some end-edge $\alpha$ of $G$ and $H-b+\beta \simeq G$ for some new end-edge $\beta$.

Therefore let us assume that this is the case and let us proceed to reconstruct, this time from $T-ab$ and $T-ac$, where $c$ is the other neighbour of $a$ in $G$. 

Reconstruction will proceed as above unless we cannot tell whether the larger component of $T-ac$ is $G$ plus an end-edge or $H$ plus an end-edge. But this ambiguity can only happen if $H+bd-\gamma \simeq G$ for some end-edge $\gamma$ of $H$ and $G-a+\delta \simeq H$ for some new end-edge $\delta$.

\begin{figure} 
 \centering
 \includegraphics[width=8.5cm, height=4cm]{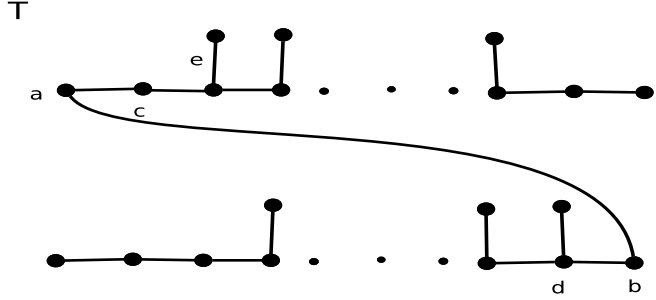}
 \caption{The tree $T$ when the two centroidal components are a conjugate pair} \label{fig:pannone0}
 \end{figure}

But this means that $G$ and $H$ are conjugate pairs and, by Theorem \ref{thm:pannone}, $T$ is therefore as shown in Figure \ref{fig:pannone0}. But then $T$ is reconstructible from $T-ab$ and $T-e$, where $e$ is as shown in Figure \ref{fig:pannone0}. Therefore $ern(T)=2$.
\lesta

\section{Edge-reconstruction number 3: three infinite families}

Molina, in \cite{molina93} had stated that $\ern(P_n) = 3$ if $T$ is a path with four or more edges. 
We shall show that his statement is correct provided that $n$, the number of vertices, is odd, that is, $P_n$ is unicentroidal.
In the following theorem we shall show that $\ern(P_n) = 2$ when $n$ is even. We shall also show that $\ern(P_n) = 3$ when $n$ is odd. To do this second part we need to show that, for each pair of cards in the edge-deck $\ED(P_n)$, there exists a graph $H\not\simeq P_n$ which has the same  pair of edge-cards in its edge-deck, that is, $H$ is a blocker for that particular pair of edge-cards.

\begin{Theorem} If  $n$ is even then $\ern(P_n) = 2$ while if $n$ is odd then $\ern(P_n) = 3$
\end{Theorem}
\noindent
{\bf Proof.}
Consider the graph $P_n$, $n$ even.  Let $e_1$ be the central edge of $P_n$ and $e_2$ any of the two edges adjacent to $e_1$. We claim that the two edge-cards $C_1 = P_n-e_1 = P_{\frac{n}{2}}\cup P_{\frac{n}{2}}$ and $C_2 = P_n - e_2 = P_{\frac{n}{2}+1}\cup P_{\frac{n}{2}-1}$ reconstruct $P_n$.  

By Molina's Lemma the graph to be reconstructed must be a tree.
Consider the missing edge of $P_n - e_1 $.  This edge can be made incident to (i) two endvertices of $P_n - e_1$; or (ii) two vertices of degree two in $P_n - e _1$; or (iii) one endvertex and one vertex of degree two. Case (i) gives $P_n$, and Case (ii) is impossible because no other edge-card of the resulting tree can be equal to the union of two paths.  Therefore we need only consider Case (iii).

Let $w$ be the vertex of degree three incident to $e_1$ after this edge is put back into $P_n - e _1$.  Then the second edge-card $C_2$ must be obtained by removing one of the other two edges incident to $w$.  But this will always give a component $P_k$ with $k>\frac{n}{2}+1$, which is a contradiction. This proves our claim.

\noindent\medskip
We now consider the odd path $P_n$ for $n = 2s + 1$ .  When two edge-cards are obtained by deleting the two edges incident to the central vertex, then a blocker would consist of the cycle $C_s$ union the path $P _{s+1}$.  The only exception is $P_5$ whose blocker in this case is the union of $C _3 $ and $P_2$.  We therefore consider any other pair of deleted edges.  Let the edges of $P_n$ be ordered as

\[e_1,e_2, \ldots, e_{n-1}.\]

Suppose we are given the two cards $P_n-e _i$ and $P _n-e _j$, where $i \leq j$. (We can assume, by symmetry, that $j\leq s$. Also, we may assume that we do not have $i=j=s$, corresponding to $i=(n-1)/2$ and $j=(n+1)/2$ since we have already observed that the blocker then is $C_s\cup P_{s+1}$.) The blocker will then consist of $S_{p,q,r}$ where $p = i$, $q=j$ and $r= 2s-i-j$. 

Therefore $\ern(P_n) > 2$. But Molina has shown that for any tree $T$ on at least four edges $\ern(T)\leq 3$, therefore $\ern(P_n) = 3 $ when $n$ is odd.
\lesta

We now show that a class of caterpillars also has $\ern=3$. 

\begin{Theorem}
The caterpillars $C(2,0,\ldots,0,2)$ of even diameter greater than 3 have edge-reconstruction number equal to 3. 
%Let $T$ be a caterpillar of order  $2n + 5$, $n\geq 1$, of the type 
%$C(2,0{^{2n - 1}},2)$. Then $\ern(T) = 3$.
\end{Theorem}

\noindent{\bf Proof.} 
Let $C$ = $C(2,0,\ldots,0,2)$ have even diameter $d>3$.
Let $(v_{0},\dots,v_{d})$ be a longest path of $C$.
By the the first result of Molina, $\ern(C)\leq 3$, so we only have to prove that $\ern(C)>2$. Thus,
we have to prove that for every pair of edge-cards $A$ and $B$ of
$C$ ($A$ and $B$ might be isomorphic), there is a blocker, that is, a graph $X$, non-isomorphic to $C$, having
two edge-cards isomorphic to edge-cards $A$ and $B$, respectively.

Let $F_{i}$ be the forest obtained by deleting edge $v_{i-1}v_{i}$,
$i=1,...,d$. Note that, because of symmetry, we need only consider
$F_{1},\ldots,F_{d/2}$. For $d>5$ we argue as follows:
\begin{itemize}
\item If the pair $F_{1}$, $F_{1}$ is chosen, we construct the graph $X$
by adding to $F_{1}$ edge $v_{0}v_{d-1}$. From $X$, we obtain $F_{1}$,
by deleting, of course, $v_{0}v_{d-1}$, and also by deleting $v_{d-1}v_{d}$.
Note that $X$ is a tree.
\item If the pair $F_{1}$, $F_{i}$ is chosen, $2\leq i\leq d/2$, we
construct graph $X$ by adding to $F_{1}$ the edge $v_{0}v_{d-i-1}$.
From $X$, we obtain  $F_{i}$
by deleting $v_{d-i}v_{d-i+1}$. Also in these cases, $X$ is a tree. 
\item If the pair $F_{j}$, $F_{i}$ is chosen, $j=2,...,(d/2)-1$,
$j\leq i\leq d/2$, we construct $X$ by adding to $F_{j}$
the edge $v_{j-1}v_{d-i}$. From $X$ we obtain $F_{i}$ by deleting 
$v_{d-i}v_{d-i+1}$. Again,
$X$ is a tree.
\item If the pair $F_{d/2}$, $F_{d/2}$ is chosen, construct $X$
by adding to $F_{d/2}$ the edge $v_{(d/2)+1}v_{d}$. From $X$ we obtain
$F_{d/2}$ by deleting the edge $v_{d-1}v_{d}$. In this last case, $X$
is not a tree, and it can be seen that there is no tree, non-isomorphic
to $C$, having $F_{d/2}$ as two of its edge-cards.
\end{itemize}
For $d=4$, we have the caterpillar $C(2,0,2)$ which we have already noted that it 
has $\ern=3$. For completeness' sake we give the same analysis as for $d>5$ above. 
\begin{itemize}
\item If the pair $F_{1}$, $F_{1}$ is chosen, we construct $X$
by adding to $F_{1}$ the edge $v_{0}v_{3}$. From $X$, we obtain $F_{1}$,
by deleting, of course, $v_{0}v_{3}$, and also by deleting $v_{3}v_{4}$.
In this case $X$ is a tree.
\item If the pair $F_{1}$, $F_{2}$ is chosen, we construct $X$ by adding
to $F_{1}$ again the edge $v_{0}v_{3}$. From $X$, we obtain 
$F_{2}$ by deleting $v_{2}v_{3}$.
The graph $X$ is a tree in this case too.
\item If the pair $F_{2}$, $F_{2}$ is chosen, we construct $X$ by adding
to $F_{2}$ the edge $v_{0}x$, where $x$ is the other edge of degree
1 in the same connected component as $v_{0}$. From $X$, we obtain
$F_{2}$ both by deleting $v_{0}x$ , and by deleting $v_{1}x$. In
this case $X$ is a not a forest.
\end{itemize}
\lesta

\medskip
[{\bf Comment.} The caterpillars $C(2,0,...,0,2)$ of odd diameter
$d$ all have $\ern=2$. Indeed, it can be directly verified that (with the same
notation as before) the pair $F_{(d-1)/2}$, $F_{(d+1)/2}$ is a pair
of edge-cards which are not in the edge-deck of any other graph not isomorphic to  $C(2,0,...,0,2)$. This observation and also Rivshin's computer search show that $C(2,0,0,0,0,2)$, which is not covered by our previous results, does indeed have $\ern=2$.]

Finally, we note that the infinite family of trees $T_k$ ($k\geq 2$) shown in Figure \ref{fig:wendyinfinitefamily} also has $\ern=3$. (Note that, when $k=2$, $T_k$ is the caterpillar $C(2,0,2)$ and, when $k=3$, $T_k$ is the graph $G_1$ shown in Figure \ref{fig:G1G2}(a).) Since there are only two types of edges up to isomorphism in $T_k$, it is easy to verify that $\ern(T_k)=3$. For example, if $e_1$ and $e_2$ are two edges of $T_k$ incident to the central vertex then $T_k-e_1$ and $T_k-e_2$ are isomorphic. The blocker having two copies of these graphs in its edge-deck is $T_{k-1} \cup R$, where $R$ is a triangle. Therefore these two subgraphs do not reconstruct $T_k$.

\section{Empirical Evidence}

\begin{figure}  
 \centering
 \includegraphics[width=8.5cm, height=4cm]{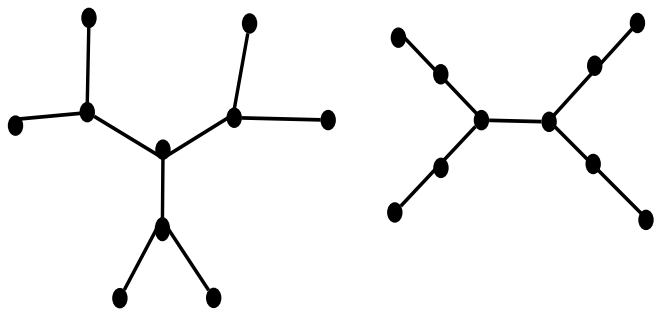}
 \caption{The graphs: (a) $G_1$; and (b) $G_2$} \label{fig:G1G2}
 \end{figure}

Empirical evidence, which was provided to us by David Rivshin \cite{rivshin11}, showed that out of more than a billion graphs on at most eleven vertices and at least four edges, only seventeen trees have edge-reconstruction number equal to 3. Four of these trees are paths of odd order which we have already considered in the previous section. Other trees are the graphs $S_{2,2,2}$, $S_{3,3,3}$ which were already noticed by Harary and Lauri \cite{harary&lau88}. Nine other trees are the caterpillars $C(2,2)$, $C(2,0,2)$, $C(1,0,1,0,1)$,  $C(2,1,2)$, $C(2,0{^3},2)$, $C(2,3,2)$, $C(2,1,1,2)$, $C(1,0,1,0,1,0,1)$ and $C(2,0{^5},2)$,  while the remaining two trees are $G_1$ and $G_2$ shown in Figure \ref{fig:G1G2}.

One can  notice that only three out of the seventeen trees are bicentroidal namely the two caterpillars $C(2,2)$ and $C(2,1,1,2)$, and the graph named $G_2$ in Figure \ref{fig:G1G2}. These trees do not contradict Theorems \ref{thm:main1} and \ref{thm:main2} since, in all three cases, the centroidal components are both pseudopaths. These small examples show that the condition that not both centroidal components are a pseudopath is required for $\ern$ to be equal to 2. Therefore the only bicentroidal trees $T$ for which we have not determined their $\ern$ because they are not covered by Theorem \ref{thm:main1} or Theorem \ref{thm:main2} are: (i) those with both centroidal components equal to pseudopaths; or (ii) those with one centroidal component being a pseudopath and the centroidal vertex in the other component having degree 2.    In the next section, we shall return to the arbitrarily large instances of these two cases, that is, when the pseudopaths involved are paths. But for the smaller cases we now present our computer search which not only covers these cases but also gives empirical evidence for our later conjecture on unicentroidal trees.

Rivshin's computer analysis considered all graphs and went up to order 11. Here, by considering only trees we extend the analysis up to order 23 (24 may be done soon, 25 and 26 are probably feasible). Because we think that it has independent interest, we shall briefly describe how this search was carried out.

\subsection{The computer search}

The program geng distributed with Brendan McKay's program
nauty \cite{mckay90} was used to generate the trees on up to 20 vertices
and nauty was used for isomorphism testing.
For trees on 21-24 vertices, Li and Ruskey's program \cite{li&rus99} was used
because it generates the trees much faster.
The approach applied to determine the edge reconstruction
number of each tree $T$ works as follows:
For each way to select two different edges $e_1$ and $e_2$ of $T$,
first create the set $S_1$ of all graphs having a card $T - e_1$ 
by adding one edge back to $T - e_1$ in all possible ways.
The number of ways to add back an edge is $n ( n-1) / 2 - (n-2)$.
Similarly, determine the set $S_2$ of graphs having a card isomorphic
to $T - e_2$.
These graphs are put into their canonical forms using nauty
(two isomorphic graphs have the same canonical form).
Next find the intersection $S$ of  $S_1$ and $S_2$ which is equal
to the set of all graphs having both cards. If the two cards
are not isomorphic to each other and $S$ 
only contains one graph then return the message that
$\ern(T)$ is equal to two. If the two cards are isomorphic to
each other, then remove from $S$ any graphs having only
one card isomorphic to $T - e_1$.
If $| S | $ is equal to one after removing these graphs then
return the message that $\ern(T)$ is equal to two. 
If all pairs of edges are tested without determining
that $ern(T)$ is equal to two, return the message
that $\ern(T)$ is equal to three.

\begin{figure}  
 \centering
 \includegraphics[width=6.5cm, height=4cm]{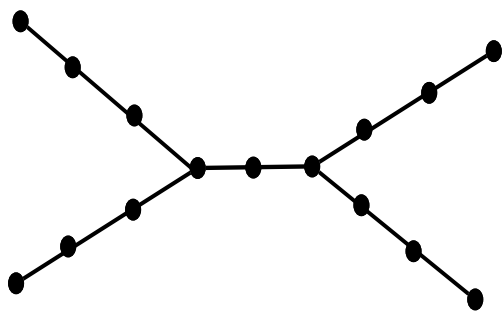}
 \caption{A tree on fifteen vertices with $\ern=3$} \label{fig:wendyG15}
 \end{figure}

The results of this computer search match the results that came from David Rivshin's data
for up to 11 vertices. This search also enabled us to discover the infinite family of trees $T_k$ with $\ern=3$ which we described above. We also found the graph $G_{15}$ on fifteen vertices (shown in Figure \ref{fig:wendyG15}) which does not fall within any known infinite class but which also has $\ern=3$.

\begin{figure}  
 \centering
 \includegraphics[width=6.5cm, height=5cm]{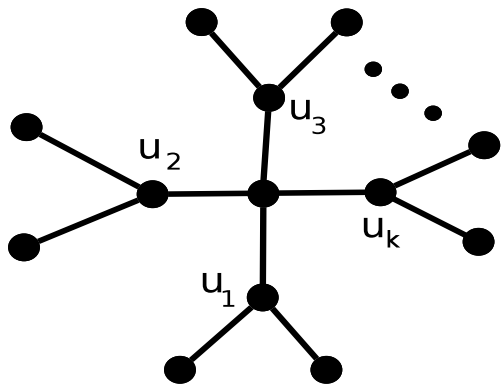}
 \caption{An infinite family of trees $T_k$ with $\ern=3$} \label{fig:wendyinfinitefamily}
 \end{figure}

\subsection{Bicentroidal trees: the remaining cases}

Let us now take stock of the situation for bicentroidal trees in the light of the results we have presented. We can summarise the situation as follows. If both centroidal components are not pseudopaths, then $\ern(T)=2$. If one component is $S_1$ then $\ern(T)=2$ except when the other component is also $S_1$ and $T$ is the caterpillar $C(2,1,1,2)$, in which case $\ern(T)=3$. If both components are $S_2$ then $\ern(T)=2$. The case when only one centroidal component is $S_2$ is covered by our computer search which confirms that, in this case too, all these trees have $\ern=2$.  

Now, if one of the two centroidal components is $P_k$, for $k\leq 5$, then $\ern(T)=2$ except when both components are $P_3$ and therefore $T$ is the caterpillar $C(2,2)$, and when the two components are $P_5$ and therefore $T$ is the graph $G_2$ of Figure \ref{fig:G1G2}. We shall therefore consider next the case when both components are $P_k$ for $k>5$.
  
First a bit of notation: Let $T$ be a bicentroidal tree with centroidal edge $ab$ and such that the two components of $T-ab$ are both isomorphic to $P_k$, the path on $k$ vertices. Let the two paths starting from $a$, but not counting $a$, have $p$ and $q$ vertices, and similarly for $b$, let the lengths be $r$ and $s$. Then we say that $T$ is of \emph{type} $T(p,q;r,s)$. We shall only give a sketch of the proof of the following lemma.

\begin{Lemma}
 Let $T$ be a bicentroidal tree with both centroidal components equal to $P_k, k>5$. Then $\ern(T)=2$.
\end{Lemma}

\noindent 
{\bf Proof.} As usual, let the centroidal edge be $ab$. It is clear that $T-ab$ cannot be one of two edge-cards giving reconstruction of $T$. Therefore let us first consider the case when we delete a non-centroidal edge $ax$ incident to $a$ and an edge $uv$ where $u$ is the endvertex on the same path from $a$ passing through $x$. Again, Lemma \ref{lem:molina} gives that $T$ is a tree. The isolated vertex $u$ must be joined to the rest of $T-uv$ in such a way that the resulting tree has an edge-card isomorphic to $T-ax$. Obviously, one way this can happen is if $u$ is joined to $v$, and this gives $T$. But there are two ``wrong'' ways in which $u$ can be joined to $T-uv$ such that the edge-card isomorphic to $T-ax$ can be obtained. Firstly, (i) $u$ can be joined to an endvertex $w$ of $T-uv$ different from $v$; or (ii) $u$ is joined to vertex $x$ of $T-uv$ . 

In case (i) let, for example, $w$ be an endvertex of the other path $P_k$ giving $T'$. Then $T'-ay$, where $y$ is the remaining neighbour of $a$, will be the edge-card isomorphic to $T-ax$; but this can happen only if $T$ is of the form $T(p,p;p,p)$. In case (ii), $T'-ax$ will be the edge-card isomorphic to $T-ax$; and this can happen only if $T$ is of the form $T(p,p;r,s)$, since $T$ is not the caterpillar $C(2,2)$ nor the graph $G_2$ of Figure \ref{fig:G1G2}.

%\emph{This last paragraph might need fixinbg. See Kevin's comments in email.}
So now we need to consider separately the case when $T$ is of the form $T(p,q;r,s)$ with $q=p$. Let $a$ be the centroidal vertex joining the two paths $P_k$, as above let $x$ be a neighbour of $a$ different from the other centroidal vertex, and let $x'$ be the other neighbour of $x$; $x'$ exists since $T$ is not the graph $G_2$ in Figure \ref{fig:G1G2}. 
We shall, in this case, use the edge cards $T-ax$ and $T-xx'$. Checking all the possibilities one finds that, again since $T$ is not the graph $G_2$,  the only way to join the two components of $T-xx'$ in such a way that the resulting tree has an edge-card isomorphic  to $T-ax$ is by joining the vertices $x$ and $x'$ in $T-xx'$.
\lesta 

\medskip\noindent
[{\bf Comment.} We think that $T$ can be reconstructed in all cases from the edge-cards $T-ax$ and $T-xx'$, as in the last paragraph of the above proof. However, reducing the problem first to the case when $q=p$ reduces the number of ways one can join the two components of $T-xx'$. This makes checking the proof shorter and easier.]

\bigskip\noindent
Therefore the only remaining case of an infinite class of bicentroidal tree whose $\ern$ is not known is when one of the centroidal components is the path $P_k$ and the other component is not a path or the graphs $S_1$, $S_2$, that is, not a pseudopath. We can now easily deal with this case in our final result which therefore neatly complements our first recosntruction result, Theorem \ref{thm:main1}.

\begin{Lemma} \label{lem:remains}
Let $T$ be a bicentroidal tree one of whose centroidal components is the path $P_k$ while the other component is not a pseudopath. Then $\ern(T)=2$.
\end{Lemma}

\noindent{\bf Proof.} As usual, let the two centroidal components of $T$ be $G$ and $H$ containing, respectively, the centroidal vertices $a$ and $b$. Suppose $H$ is the path $P_k$. Therefore, since $G$ is not a pseudopath, we may assume that $\deg(b)=2$, otherwise we know that $\ern(T)=2$, by Theorem \ref{thm:main1}. Let $d$ be the other neighbour of $b$, therefore $d\in V(H)$. We shall first try to reconstruct $T$ from $T-ab$ and $T-bd$.  

As usual, we can determine from these two edge-cards that $T$ must be a tree, by Lemma \ref{lem:molina}, and we therefore know the two centroidal components of $T$. 
We now consider $T-bd$. The two usual situations can arise: (1) the large component of $T-bd$ is isomorphic to $G$ plus an edge and the smaller component is isomorphic to $P_k$ less an edge; or (2) vice-versa, with the roles of $G$ and $P_k$ reversed. Suppose the first but not the second case holds. Therefore we need to find, in the large component $K$ of $T-bd$, an end-edge $e$ such that $K-e$ is isomorphic to $G$. Reconstruction would then proceed by joining the endvertex incident to $e$ with an endvertex of the other component of $T-bd$. The edge $ab$ is surely one such edge $e$. But even if there is another edge $vw$ with $\deg(v)=1$ such that $K-vw\simeq G$, then $w$ and $b$ would be similar in $K$. Therefore, joining the vertex $w$ or the vertex $b$ to an endvertex of the other component would give isomorphic reconstructions. This shows that $T$ is reconstructible from $T-ab$ and $T-bd$.

We can now suppose that case (2) holds. This means that $G+ab - \alpha \simeq P_k$ for some end-edge $\alpha$ of $G$.
So, since $G$ is not $P_k$, it must be the path $P_{2k-1}$ with an extra end-edge incident to one of its interior vertices. Therefore $T$ is the path $P_{2k-1}$ plus an end-edge incident to some interior vertex. 

But in this case it is easy to shown that $\ern(T)=2$. Let the vertices of the path $P_{2k-1}$
be, in order, $v_1,v_2,\ldots,v_{2k-1}$. Let the extra end-edge be $v_ix$ where $i$ is not equal to 1 or $2k-1$. Since $T$ is bicentroidal, $v_i$ is not the central vertex $v_k$ of $P_{2k-1}$. Also, we may assume, without loss of generality, that $i<k$. But then it is easily checked that $T$ is reconstructible from the edge-cards $T-v_ix$ and $T-v_{2i-1}v_{2i}$.
\lesta

This final result and the previous comments gives the main result of this paper.

\begin{Theorem} \label{thm:main3}
Every bicentroidal tree except $C(2,2)$, $C(2,1,1,2)$ and the graph $G_2$ shown in Figure \ref{fig:G1G2}(b) has edge-reconstruction number equal to 2.
\end{Theorem}

\section{Final comments}

We have managed to fill in the gaps in our knowledge of the edge-recons\-truction number of bicentroidal trees. The computer search described above also leads us to make this conjecture for unicentroidal trees.

\begin{Conjecture}
The only infinite classes of trees which have $\ern=3$ are the paths on an odd number of vertices, the caterpillars $C(2,0,\ldots,0,2)$ of even diameter, and the family of trees $T_k$ depicted in Figure \ref{fig:wendyinfinitefamily}.
\end{Conjecture}

Proving this conjecture might not be easy. The difficulty of determining which unicentroidal trees have $\ern$ equal to 2 or 3 when the (vertex) reconstruction number of trees is known is again evidence for the phenomenon, commented upon in \cite{asciak&10},  when determining the edge-reconstruction number of a class of graphs is sometimes more difficult than determining the (vertex) reconstruction number.

\section*{Acknowledgement}
We are again grateful to David Rivshin whose data and computer programs helped us when we could not see the trees for the wood when considering the exceptional cases amongst smaller trees.

\bibliography{ern(T)3.bib}

\end{document}